\documentclass[11pt, letterpaper]{article}

\usepackage[
letterpaper,
top=1in,
bottom=1in,
left=1in,
right=1in]{geometry}

\usepackage{pifont}
\usepackage[utf8]{inputenc}

\usepackage{amsmath, amsthm, amssymb}
\usepackage{mathtools}

\usepackage[dvipsnames, table]{xcolor}
\usepackage[T1]{fontenc}
\usepackage[full]{textcomp}
\usepackage[american]{babel}

\usepackage{caption}

\usepackage{array,booktabs,longtable}
\usepackage{colortbl}

\definecolor{SOSstripe}{gray}{0.94}
\newcolumntype{Z}[1]{>{\raggedright\arraybackslash}m{#1}}
\newlength{\SOSwidth}

\usepackage{siunitx}

\usepackage[pagebackref,colorlinks=true,urlcolor=blue,linkcolor=blue,citecolor=OliveGreen]{hyperref}
\usepackage[nameinlink]{cleveref}
\crefname{equation}{}{}

\usepackage{enumitem}
\usepackage{multirow}
\usepackage{graphicx}

\newtheorem{theorem}{Theorem}[section]

\newtheorem{proposition}[theorem]{Proposition}

\theoremstyle{definition}
\newtheorem{definition}[theorem]{Definition}

\newtheorem{example}[theorem]{Example}
\newtheorem{question}[theorem]{Question}
\newtheorem{remark}[theorem]{Remark}

\Crefname{question}{Question}{Questions}

\newtheorem*{question*}{Question}

\newcommand{\Lo}{\text{\L}}
\newcommand{\R}{\mathbb{R}}
\newcommand{\N}{\mathbb{N}}

\newcommand{\x}{\mathbf{x}}
\newcommand{\mg}{\mathbf{g}}

\newcommand{\preo}{\mathcal{T}}
\newcommand{\qm}{\mathcal{Q}}
\newcommand{\ideal}{\mathcal{I}}
\newcommand{\pos}{\mathcal{P}}
\newcommand{\se}{\mathcal{S}}

\newcommand{\M}{\mathrm{Mtz}}
\newcommand{\CL}{\mathrm{CL}}
\newcommand{\BS}{\mathrm{BS}}

\newcommand{\e}{\varepsilon}
\newcommand{\eps}{\varepsilon}

\newcommand{\gcube}{\mg_{\rm cube}}

\newcommand{\UPPER}{upper bounds on degree and error}
\newcommand{\LOWER}{lower bounds on degree and error}

\newcommand{\lean}{\textbf{\textcolor{OliveGreen}{Lean}}}

\usepackage{authblk}

\makeatletter
\newcommand{\AB@affilfootnotes}{}
\DeclareRobustCommand{\AB@printaffilfootnotes}{\AB@affilfootnotes}
\g@addto@macro\@thanks{\AB@printaffilfootnotes}

\renewcommand{\affil}[2][]{%
  \g@addto@macro\AB@affilfootnotes{%
    \begingroup
      \renewcommand{\thefootnote}{#1}%
      \footnotetext[0]{#2}%
    \endgroup
  }%
}

\renewcommand{\@author}{%
  \def\rlap##1{##1}
  \AB@authlist
}
\makeatother

\title{Optimal degree bounds for sum-of-squares representations}

\author[1]{Alexander Taveira Blomenhofer}

\author[2,*]{Sander Gribling}

\author[3,4]{Didier Henrion}

\author[3,4]{Milan Korda}

\author[5]{Poornendu Kumar}

\author[3]{Victor Magron}

\author[6,*]{Lucas Slot}

\author[3]{Luis~Felipe Vargas}

\author[2]{Juan Vera}

\affil[1]{University of Porto}
\affil[2]{Tilburg University}
\affil[3]{LAAS CNRS}
\affil[4]{Czech Technical University in Prague}
\affil[5]{Faculty of Mathematics and Physics, University of Ljubljana}
\affil[6]{University of Amsterdam.}
\affil[*]{Corresponding authors. Emails: \url{s.j.gribling@tilburguniversity.edu}, \url{l.f.h.slot@uva.nl}}

\begin{document}
\maketitle
\vspace{-1cm}
\begin{abstract}
The moment-sum-of-squares (SOS) hierarchy provides a sequence of bounds on the minimum of a polynomial over a basic closed semialgebraic set, i.e., a subset of $\R^n$ defined by finitely many polynomial constraints.
We establish optimal bounds on the error of the SOS hierarchy in the following settings:
\begin{enumerate}
    \item If the constraints define an Archimedean quadratic module: the error of the Putinar-based hierarchy at level $k$ is $\Theta((\log k)^4/k^2)$ in the worst case.
    \item If the constraints define a compact semialgebraic set: the error of the Schmüdgen-based hierarchy at level $k$ is $\Theta((\log k)^2/k^2)$ in the worst case.
    \item If the constraints define an Archimedean quadratic module and satisfy linear independence constraint qualification conditions (LICQ): the error of the hierarchy at level $k$ is~$\Theta(1/k^{4})$ in the worst case. Here, the upper bound holds already for the Putinar-based hierarchy, while the lower bound holds even against the Schmüdgen-based hierarchy.
\end{enumerate}
Equivalently, these give best-possible effective versions of Putinar's and Schm\"udgen's Positivstellens\"atze for these settings (in terms of their dependence on $k$).
\end{abstract}

\setcounter{tocdepth}{2}
\tableofcontents

\section*{Statement}
This is an announcement of results which were produced over the course of several months with extensive use of generative AI. Several subsets of the authors independently arrived at (parts of) these results over that time period either by prompting various models, by independent human thought or, in most cases, by a mixture of both. 
The purpose of this announcement is to have a common, comprehensive reference for a set of results which we believe to essentially settle a line of research going back more than 20 years. 

For each of the results claimed below, we provide a certificate of correctness in Lean. These can be viewed at \url{https://github.com/Lucas-Slot/OptimalConvergenceRates}. All theorem statements marked with {\lean}~have a direct formalization in that repository. This announcement does not contain proofs beyond these certificates.
In particular, we have opted \emph{not} to publish any proofs generated directly by AI. Our motivation for omitting such proofs is twofold: 
\begin{itemize}
    \item They are not suited towards furthering human understanding of the underlying math;
    \item They do not properly incorporate, reference, or credit existing ideas from the literature.
\end{itemize}

The authors are preparing expositions of all claimed results that take into account the two points above. We intend to provide these in a future version of this preprint.
\emph{We would like to do so as a community, and we welcome contributions from others who may have been working on these or similar problems}.

\section{Introduction}
Let $\R[\x]$ denote the polynomial ring in $n$ variables $\x = (x_1, x_2, \ldots, x_n)$. For a tuple of polynomials~$\mg = (g_1, g_2, \ldots, g_m)$, we write $\se(\mg)$ for the (basic, closed) semialgebraic set 
\[
    \se(\mg) := \{ \x \in \R^n : g_1(\x) \geq 0, \ldots, g_m(\x) \geq 0 \}.
\]
We consider the classical problem of certifying whether a given polynomial $p \in \R[\x]$ belongs to the cone $\pos(\mg)$ of nonnegative polynomials on $\se(\mg)$. 
That is, whether $p(\x) \geq 0$ for all $\x \in \se(\mg)$.
One way to certify \emph{global} nonnegativity of $p$ is to write it as a \emph{sum of squares}, i.e.,
\[
    p(\x) = q_1(\x)^2 + \ldots + q_\ell(\x)^2 \quad (q_i \in \R[\x]).
\]
We write $\Sigma[\x] \subseteq \R[\x]$ for the cone of all polynomials admitting such a decomposition.
This idea can be extended to certify nonnegativity over $\se(\mg)$ as follows. For $k \in \N$, consider the (truncated) quadratic module and preordering associated to $\mg$, given by
\begin{align}
    \qm(\mg)_{k} &:= \left\{\sum_{i=0}^m g_i\sigma_i \,:\, \sigma_i \in \Sigma[\x] \text{ and } \mathrm{deg}(g_i \sigma_i) \leq 2k \text{ for all } 0\leq i \leq m \right\}, \\
    \preo(\mg)_{k} &:= \left\{ \sum_{I \subseteq [m]} g_I \sigma_I \,:\, \sigma_I \in \Sigma[\x] \text{ and } \mathrm{deg}(g_I \sigma_I) \leq 2k \text{ for all } I \subseteq [m] \right\}.
\end{align}
Here, $g_I := \prod_{i \in I} g_i$ for $I \subseteq [m] = \{1,2,\ldots,m\}$, and we have adopted the notational convention that ${g_0 = g_\emptyset = 1}$. We write $\qm(\mg) =\bigcup_{k=0}^\infty \qm(\mg)_k$ and  $\preo(\mg) = \bigcup_{k=0}^\infty \preo(\mg)_k$ for the degree-unrestricted versions, respectively.
It follows immediately from the definition that
\begin{equation} \label{EQ:preoquadinc}
    \qm(\mg)_k \subseteq \preo(\mg)_k \subseteq 
    \pos(\mg) \quad (\forall k\in \N).
\end{equation}
Krivine \cite{krivine_anneaux_1964} and Stengle
\cite{stengle_nullstellensatz_1974}
showed that any nonnegative polynomial on $\se(\mg)$ can be represented as a \emph{ratio} of polynomials in the preordering. This generalizes Artin's solution to Hilbert's seventeenth problem \cite{artin_uber_1927}.
A natural question is then whether all nonnegative polynomials on~$\se(\mg)$ lie in $\preo(\mg)$ or even in $\qm(\mg)$. This turns out not to be true in general. However, the classical \emph{Positivstellens\"atze} of Putinar and Schm\"udgen tell us that the inclusions~\eqref{EQ:preoquadinc} are ``asymptotically tight'' for (algebraically) compact (basic) semialgebraic sets, in the following sense:
\begin{theorem}[{Putinar's Positivstellensatz~\cite{Putinar1993}}]\label{THM:putinar}
Suppose that $\qm(\mg)$ is Archimedean.\footnote{Meaning that $\qm(\mg)$ contains the ball constraint $N - \|\x\|^2$ for some $N > 0$. In particular, $\se(\mg)$ is then compact.} Let~${p \in \pos(\mg)}$. For every $\e > 0$, there exists $k \in \N$ such that
$
    p + \e \in \qm(\mg)_k.
$
\end{theorem}
\begin{theorem}[{Schm\"udgen's Positivstellensatz~\cite{Schmudgen1991}}]\label{THM:schmudgen}
Suppose that $\se(\mg) \subseteq \R^n$ is compact. Let~$p \in \pos(\mg)$. For every $\e > 0$, there exists $k \in \N$ such that
$
    p + \e \in \preo(\mg)_k.
$    
\end{theorem}
In particular, these imply that any polynomial \emph{strictly} positive on $\se(\mg)$ lies in $\qm(\mg)$ and $\preo(\mg)$. 
\subsection{Effective Positivstellens\"atze and convergence of the SOS hierarchy.}

There has been substantial interest in proving 
(asymptotic) bounds on the degree $k = k(\e)$ required in the above theorems as $\e \to 0$. This is motivated in large part by the fact that such bounds translate directly to convergence rates for the celebrated moment-SOS hierarchy for polynomial optimization~\cite{lasserre2001global,parrilo2000thesis}. This can be seen as follows:
For a nonempty compact feasible set $\se(\mg)$ and an objective
$f\in\R[\x]$, put $p=f-f_{\min}$, where $f_{\min} := \min_{\x \in \se(\mg)} f(\x)$. Note that $p \in \pos(\mg)$ and has zeros in $\se(\mg)$.
For
$\mathcal C\in\{\qm(\mg),\preo(\mg)\}$, define
\begin{align}
f_{k,\mathcal C}
&:= \sup\{\lambda\in\R:f-\lambda\in\mathcal C_k\},\\
\e_{p,k,\mathcal C}
&:= \inf\{\e\ge0:p+\e\in\mathcal C_k\}
= f_{\min}-f_{k,\mathcal C}.
\label{eq:epkC}
\end{align}
The parameter $\e_{p,k,\mathcal C}$ is called the \emph{error} of the moment-SOS hierarchy at level $k$. Its asymptotic behavior as $k \to \infty$ is referred to as the \emph{convergence rate} of the hierarchy.
Equivalently, in the language of~\Cref{THM:putinar,THM:schmudgen}, $\e_{p,k,\mathcal C}$ corresponds to a degree requirement 
\begin{equation}
 k_{p,\e,\mathcal C}
 :=\min\{k\in\N:p+\e\in\mathcal C_k\}.
\end{equation}
The central question we address in this paper is the following:
\begin{question} \label{question}
Under suitable assumptions on $\mg$, what is the dependence of $\e_{p,k,\mathcal C}$ on $k$ as $k \to \infty$ (or, equivalently, of $k_{p,\e,\mathcal C}$ on $\e$ as $\e \to 0$)?
\end{question}
This question has been well-studied in the literature.
Below, we give a brief overview of the current state of the art. We end the overview with \Cref{tab:SOTA}, which compares the state of the art to the new results provided by this work. We present our results in detail in \Cref{sec:contributions}. Finally, we give a (hopefully) complete history of the development in~\Cref{sec:tables} for reference.

\subsection{State of the art}
\subsubsection{Upper degree bounds and convergence rates.}
Upper bounds on the error $\e_{p,k,\mathcal C}$, 
or equivalently, on the required degree $k_{p,\e,\mathcal C}$
provide a safety guarantee that
the hierarchy cannot converge slower than a stated rate, even on the hardest instances in the class.

The literature on these upper bounds can be roughly divided into two streams. On the one hand, there is a line of works establishing such bounds under general constraint sets $\mg$. 
Roughly speaking, this has culminated in the following: For every $\mg$ defining an Archimedean quadratic module, there is a constant $c(\mg) > 0$, so that, for every $p \in \pos(\mg)$,
\[
\e_{p,k,\qm(\mg)} \leq O(  k^{-1/c(\mg)}).
\]
Here, the $O$-notation hides constants depending only on $p, \mg$.\footnote{In many of the works discussed below, these constants are made (more or less) explicit. For simplicity, we do not specify them in this overview.}
The same is true for $\preo(\mg)$ under the weaker assumption that $\se(\mg)$ is compact (but with a potentially larger constant $c(\mg)$). 
Thus, the dependence of $\e_{p,k,\mathcal C}$ on $k$ is at most inverse polynomial.

A long line of work has led to the following understanding of this inverse-polynomial dependence: The constant $c(\mg)$ can be bounded in terms of the so-called {\L}ojasiewicz exponent $\Lo(\mg)$ of $\mg$ (\Cref{DEF:Loja}). For Putinar ($\mathcal{C} = \qm(\mg)$), one has $c(\mg) \leq 2 \cdot \Lo(\mg)$, and for Schm\"udgen ($\mathcal{C} = \preo(\mg)$), one has $c(\mg) \leq \Lo(\mg)$~\cite{heijmans2026degree, tran2026convergence}. In the best case (when $\mg$ satisfies a type of constraint qualification, see~\Cref{def::LICQ}), we have $\Lo(\mg) = 1$, leading to $c(\mg) = 2$ and $c(\mg) = 1$, respectively.

On the other hand, there is a series of works investigating error bounds for specific sets of constraints $\mg$, corresponding to fundamental domains of optimization. The most well-studied examples are 
\begin{align}
    \mg_{\rm sph} &= (1-\|\x\|^2, \, \|\x\|^2 - 1), \text{ defining the unit sphere } \se(\mg) = S^{n-1}; \\
    \gcube &= (1 - x_i^2 : 1 \leq i \leq n), \text{ defining the unit hypercube } \se(\mg) = [-1,1]^n.
\end{align}
For these examples, one can prove much better degree bounds than the general results above would predict. 
In particular, it is known that $c(\gcube) = c(\mg_{\rm sph}) = 1/2$ 
suffice for certificates in the preordering~\cite{fang2021sum, slot2022sum}. (Note that, for the standard description of the unit sphere, the quadratic module and preordering coincide.) 
For representations in the quadratic module of the hypercube, one can select  $c(\gcube) = 1/2$, 
up to polylogarithmic factors. 

\subsubsection{Lower degree bounds and convergence obstructions.}
Much less is known about \emph{lower bounds} on $\e_{p,k,\mathcal C}$ in terms of $k$ (or, equivalently, on $k_{p,\eps,\mathcal C}$ in terms of $\e$).
Such lower bounds represent necessity statements: they quantify intrinsic limitations of
the method (or of the representation class), and rule out uniformly better guarantees without changing assumptions or algorithms. Essentially, only two lower bound constructions are known.

First, a classical example due to Stengle~\cite
{Stengle1996}
shows that degree $k = \Omega(\e^{-1/2})$ is needed to represent $1-x^2+\e$ in the preordering generated by~$(1-x^2)^3$. That is, $c(\mg) \geq 1/2$ for that example.
Recently, Henrion~\cite{Henrion2026} showed that this analysis is tight, and derived an exact relaxation order.
Stengle's construction relies fundamentally on the non-standard choice of representation for the interval.
Indeed, any polynomial $p$ nonnegative on $[-1,1]$ admits a certificate in~${\preo(1-x^2)}$ of degree at most
$\mathrm{deg}(p)+1$ (by the Markov--Luk\'acs theorem). The key point is that the constraint $(1-x^2)^3$ is degenerate at the boundary of $[-1,1]$ (cf.~\Cref{def::LICQ} below). 

Second, for semialgebraic sets with nondegenerate description, to our knowledge, the only available lower bound is due to~\cite{BaldiSlot2024}. There, the authors show that $(1-x_1^2)(1-x_2^2) + \e$ is not an element of the quadratic module~${\qm(1-x_1^2,1-x_2^2)_k}$ unless~$k = \Omega(\e^{-1/8})$. No lower bounds are available for the preordering of the unit hypercube, or for the quadratic module of the unit sphere (with standard description).

\subsection{Comparison to new results}
In \Cref{tab:SOTA} below, we give an overview of the state-of-the-art results discussed above, and compare them to the new (optimal) results obtained in this work.

\begin{table}[ht!]
    \centering
    \small
    \renewcommand{\arraystretch}{1.4}
    \setlength{\tabcolsep}{4pt}

    \begin{tabular}{@{}llccc@{}}
        \toprule
        \rowcolor{white}
        & & \multicolumn{2}{c}{Previous work}
        & \textbf{This work} \\
        \cmidrule(lr){3-4}
        Setting & Certificate
        & Upper bound $O(\cdot)$ & Lower bound $\Omega(\cdot)$
        & \textbf{Optimal rate} $\Theta(\cdot)$ \\
        \midrule
        \rowcolor{SOSstripe}
        Compact
        & $\preo(\mg)$
        & $k^{-1/\Lo(\mg)}$
        & $k^{-2}$
        & $(\log k)^2/k^2$ \\

        Archimedean
        & $\qm(\mg)$
        & $k^{-1/(2\Lo(\mg))}$
        & $k^{-2}$
        & ${(\log k)^4/k^2}$ \\

        \addlinespace[4pt]
        \rowcolor{SOSstripe}
        Hypercube
        & $\preo(\gcube)$
        & $k^{-2}$
        & ---
        & $k^{-4}$ \\

        Hypercube
        & $\qm(\gcube)$
        & $(\log k)^3/k^2$
        & $k^{-8}$
        & $k^{-4}$ \\

        \addlinespace[4pt]
        \rowcolor{SOSstripe}
        Archimedean + $\Lo(\mg)=1$
        & $\qm(\mg)$
        & $k^{-1/2}$
        & $k^{-8}$
        & $k^{-2}$ \\

        Archimedean + LICQ
        & $\qm(\mg)$
        & $k^{-1/2}$
        & $k^{-8}$
        & $k^{-4}$ \\

        \bottomrule
    \end{tabular}

    \caption{Comparison of prior upper and lower bounds on the asymptotic
    error $\e_{p,k,\mathcal C}$, as $k \to \infty$ with the optimal
    rates established in this work (Putinar = $\qm$, Schm\"udgen = $\preo$). Upper bounds hold for all
    admissible instances, whereas lower bounds are attained by
    specific instances. The rates in the last column are optimal
    in the worst-case sense up to multiplicative constants (that may depend on the dimension $n$ as well as on $p$ and $\mg$).}
    \label{tab:SOTA}
\end{table}

\section{Contributions} \label{sec:contributions}

We provide a series of matching upper and lower degree bounds for Putinar's and Schm\"udgen's Positivstellensatz in three regimes: (1) general constraint sets $\mg$ defining an Archimedean quadratic module or compact semialgebraic set, respectively; (2) the standard description $\gcube$ of the unit hypercube~$[-1,1]^n$; (3) general $\mg$ satisfying the following constraint qualification condition.

\begin{definition}\label{def::LICQ}
We say that a tuple of polynomials $\mg$ satisfies Linear Independence Constraint Qualification (LICQ) if, for every $\x \in \se(\mg)$, the gradients of the active constraints at $\x$, being
\begin{align}\label{eq:LICQ-gradients}
    {\big(\nabla g(\x) : g \in \mg,~g(\x) = 0\big)},
\end{align}
are linearly independent (in particular, nonzero). By convention, if the tuple $\mg$ contains both $g$ and~$-g$, then only one of their gradients is considered in eq. \Cref{eq:LICQ-gradients}.
\end{definition}

Notable examples of tuples that satisfy the LICQ include $\mg_{\rm sph} = (1-\|\x\|^2,\|\x\|^2-1)$ corresponding to the unit sphere, and $\mg_{\rm ball} = (1-\|\x\|^2)$ corresponding to the unit ball. 

\subsection{Effective Putinar: Archimedean quadratic modules} \label{sec:effective_putinar}
The following theorems provide an asymptotically optimal effective version of Putinar's Positivstellensatz (\Cref{THM:putinar}). Namely, they establish that degree $k = \Theta(\log(1/\e)^2 / \sqrt \e)$ is sufficient in general, and necessary in the worst case. We first show the upper degree bound. 

\begin{theorem}[{\lean} | Archimedean, \UPPER] \label{THM:generalputinarupper} Let $\mg$ be a set of constraints defining an Archimedean quadratic module $\qm(\mg)$, and a nonempty semialgebraic set $\se(\mg)$. Let $p \in \pos(\mg)$.
Then, there exist constants $C, k_0, \e_0 > 0$, depending only on $p$ and $\mg$, such that
\[
    \e_{p, k, \qm(\mg)} \leq C \cdot \frac{(\log k)^4}{k^2} \quad (\forall k \geq k_0); \quad 
    k_{p, \e, \qm(\mg)} \leq C \cdot \frac{\log(1/\e)^2}{\sqrt{\e}} \quad (\forall \e \leq \e_0).
\]

\end{theorem}

The lower degree bound uses the following instance. 
\begin{example} \label{EXMP:generalputinar}
Let $n = 3$ and consider the constraint set
\[
\mg=(x_1,x_2, \, x_1x_2x_3+x_3^5, \, 1-\|\x\|^2).
\]
The quadratic module $\qm(\mg)$ is Archimedean. The semialgebraic set is $\se(\mg) = \{ \x:\x \geq 0, \, \|\x\|^2 \leq 1 \}$, i.e., the intersection of the unit ball with the nonnegative orthant. To see this, note that given the first two constraints, the third is equivalent to $x_3 \geq 0$. In particular, $\se(\mg)$ is nonempty, and thus $\mg$ satisfies the conditions of~\Cref{THM:generalputinarupper}.
\end{example}

\begin{theorem}[{\lean} | Archimedean, \LOWER]  \label{THM:generalputinarlower} Let $n=3$, and let $\mg$ be as in~\Cref{EXMP:generalputinar}. Consider the polynomial $p = x_3$, which is nonnegative on $\se(\mg)$. There exist constants $c, k_0, \e_0 > 0$ such that
\begin{equation}
 \e_{p, k, \qm(\mg)} \geq c \cdot \frac{(\log k)^4}{k^2} \quad (\forall k \geq k_0); \quad  k_{p, \e, \qm(\mg)} \geq c \cdot \frac{\log(1/\e)^2}{\sqrt{\e}} \quad (\forall \e \leq \e_0);
\end{equation}
\end{theorem}

\begin{remark}
    The constraint set $\mg$ of~\Cref{EXMP:generalputinar} does not satisfy LICQ. Indeed, for $g_3(x_1,x_2,x_3) = x_1x_2x_3 + x_3^5$, we have  
        $g_3(0,0,0) = 0$ and $\nabla g_3(0,0,0) = (0,0,0)$.
\end{remark}

\subsection{Effective Schm\"udgen: Compact semialgebraic sets}
The following theorems provide an asymptotically optimal effective version of Schm\"udgen's Positivstellensatz (\Cref{THM:schmudgen}). Namely, they establish that degree $k = \Theta(\log(1/\e) / \sqrt \e)$ is sufficient in general, and necessary in the worst case. (This result differs from that of the previous section by a factor $\log(1/\e)$).

\begin{theorem}[{\lean} | Compact, \UPPER] \label{THM:generalschmudgenupper} Let $\mg$ be a set of constraints defining a compact, nonempty semialgebraic set $\se(\mg)$. Let $p \in \pos(\mg)$.
Then, there exist constants $C, k_0, \e_0 > 0$, depending only on $p$ and $\mg$, such that
\[
    \e_{p, k, \preo(\mg)} \leq C \cdot \frac{(\log k)^2}{k^2} \quad (\forall k \geq k_0); \quad 
    k_{p, \e, \preo(\mg)} \leq C \cdot \frac{\log(1/\e)}{\sqrt{\e}} \quad (\forall \e \leq \e_0).
\]
\end{theorem}
The lower bound is based on the following bivariate example. 

\begin{example} \label{EXMP:generalschmudgen}
Let $n=2$ and consider the constraint set
\[
    \mg = (x_2(x_1^2+x_2^2),\, 1-x_1^2-x_2^2).
\]
The semialgebraic set $\se(\mg) = \{x_2 \geq 0, \, \|\x\|^2 \leq 1\}$  is compact. In fact, $\qm(\mg)$ is Archimedean. 
\end{example}

\begin{theorem}[{\lean} | Compact, \LOWER]
Let $n=2$, and let $\mg$ as in~\Cref{EXMP:generalschmudgen}. Consider the polynomial $p = x_2$, which is nonnegative on $\se(\mg)$. There exist constants $c, k_0, \e_0 > 0$ such that
\[
    \e_{p, k, \preo(\mg)} \geq c \cdot \frac{(\log k)^2}{k^2} \quad (\forall k \geq k_0); \quad 
    k_{p, \e, \preo(\mg)} \geq c \cdot \frac{\log(1/\e)}{\sqrt{\e}} \quad (\forall \e \leq \e_0).
\]
\end{theorem}

\begin{remark}
    The constraint set $\mg$ of~\Cref{EXMP:generalschmudgen} does not satisfy the LICQ (\Cref{def::LICQ}). Indeed, for $g_1(x_1,x_2) = x_2(x_1^2+x_2^2)$, we have  
        $g_1(0,0) = 0$ and $\nabla g_1(0,0) = (0,0)$.
\end{remark}

\subsection{The unit hypercube}
In this section, we consider the unit hypercube $[-1,1]^n$, with its standard description
\[
    [-1,1]^n = \se(\gcube), \quad \text{ where } \gcube = (1-x_i^2 \, : \, 1 \leq i \leq n).
\]
The following theorems provide asymptotically optimal effective versions of both Putinar's and Schm\"udgen's Positivstellens\"atze as applied to the unit hypercube (with its standard description). Namely, they establish for both that degree $k = \Theta(1 / \sqrt[4] \e)$ is sufficient in general, and necessary in the worst case.
For Putinar, the worst-case example uses $n = 2$ variables. For Schm\"udgen, the worst-case example uses $n = 4$ variables.

\begin{theorem}[{\lean} | Hypercube, \UPPER] \label{THM:cubeupper}
For every $n\ge1$ and $d\ge0$, there exist constants $C_{n,d}, 
k_{n,d}, \e_{n,d} > 0$, depending only on $n$ and $d$, such that, for every polynomial $p \in \pos(\gcube)$ of degree at most~$d$,
\[
\e_{p, k, \qm(\gcube)} \leq C_{n,d} \cdot \frac{\|p\|_\infty}{k^4} \quad (\forall k \geq k_{n,d}); 
\quad k_{p, \e, \qm(\gcube)} \leq C_{n,d} \cdot \left(\sqrt[4]{\frac{\|p\|_\infty}{\e}} + 1 \right)\quad (\forall \e \leq \e_{n,d}).
\]
Here, $\|p\|_\infty := \max_{\x \in [-1,1]^n} |p(\x)|$. {As a consequence, the same holds for $\preo(\gcube)$.}
\end{theorem}

We first present the instance used in the lower degree bound for Putinar-type certificates, for $n=2$.

\begin{example} \label{EXMP:BS}
    Let $n = 2$. Consider the polynomial 
    \[
        \BS(x_1, x_2) := (1-x_1^2)(1-x_2^2).
    \]
    This polynomial is nonnegative on $[-1, 1]^2$.
    It is an element of $\preo(\gcube)$ by definition,
    but it is not an element of~$\qm(\gcube)$, see~\cite{BaldiSlot2024}.
\end{example}
\begin{theorem}[{\lean} | Hypercube, \LOWER~I] \label{THM:cubelowerI}
Let $n=2$, and let $p = \BS \in \R[\x]$ be as in~\Cref{EXMP:BS}. There exist  absolute constants $c, k_0, \e_0 > 0$ such that, 
\[
\e_{p, k, \qm(\gcube)} \geq c \cdot \frac{1}{k^4} \quad (\forall k \geq k_0); 
\quad k_{p, \e, \qm(\gcube)} \geq c \cdot  \frac{1}{\sqrt[4]{\e}} \quad (\forall \e \leq \e_0).
\]
These bounds are asymptotically tight in light of~\Cref{THM:cubeupper}.
\end{theorem}

\begin{remark}
    The statement of \Cref{THM:cubelowerI} was conjectured based on experiments with a high-precision SDP solver~\cite{BaldiKocvara:PC}. In fact, Baldi and Kočvara  conjecture it holds more generally for the representation of $g_1 \cdot g_2 + \eps$ in $\qm(g_1, g_2)$, where $g_1, g_2$ are arbitrary polynomials that intersect transversally in $\R^n$ and generate an Archimedean quadratic module. 
\end{remark}

We now present two instances that lead to lower degree bounds for Schmüdgen-type certificates for $n=3$, and $n \geq 4$.

\begin{example} \label{EXMP:Mtz}
    Let $n = 3$. Consider the Motzkin form 
    \[
        \M(x_1, x_2, x_3) := x_1^4x_2^2+x_1^2x_2^4+x_3^6-3x_1^2x_2^2x_3^2.
    \]
    The polynomial $\M$ is globally nonnegative, but it is not a sum of squares, see~\cite{Reznick2007}. In fact, it is not an element of $\preo(\gcube)$.
\end{example}

\begin{theorem}[{\lean} |  Hypercube, \LOWER~II] 
Let $n=3$, and let $p=\M \in \R[\x]$ be as in~\Cref{EXMP:Mtz}. There exist absolute constants $C, c, k_0, \e_0 > 0$ such that
\[
C \cdot \frac{1}{k^6} \geq \e_{p, k, \preo(\gcube)} \geq c \cdot \frac{1}{k^6} \quad (\forall k \geq k_0); 
\quad C \cdot \frac{1}{\sqrt[6]{\e}} \geq k_{p, \e, \preo(\gcube)} \geq c \cdot \frac{1}{\sqrt[6]{\e}} \quad (\forall \e \leq \e_0).
\]
\end{theorem}

\begin{example} \label{EXMP:CL}
    Let $n = 4$. Consider the Choi--Lam form 
    \[
        \CL(x_1,x_2,x_3,x_4)  := x_1^4+x_2^2x_3^2+x_3^2x_4^2+x_4^2x_2^2
  -4x_1x_2x_3x_4.
    \]
    The polynomial $\CL$ is globally nonnegative, but is not a sum of squares, see~\cite{Reznick2007}. In fact, it is not an element of $\preo(\gcube)$.
\end{example}

\begin{theorem}[{\lean} |  Hypercube, \LOWER~III] 
Let $n=4$, and let $p=\CL \in \R[\x]$ be as in~\Cref{EXMP:CL}. There exist  absolute constants $c, k_0, \e_0 > 0$ such that, 
\[
\e_{p, k, \preo(\gcube)} \geq c \cdot \frac{1}{k^4} \quad (\forall k \geq k_0); 
\quad k_{p, \e, \preo(\gcube)} \geq c \cdot  \frac{1}{\sqrt[4]{\e}} \quad (\forall \e \leq \e_0).
\]
Consequently, the same holds for $\qm(\gcube)$. These bounds are asymptotically tight in light of~\Cref{THM:cubeupper}.
\end{theorem}

\begin{remark}
    {We presented three different lower bounds, respectively for $n=2$, $n=3$, $n=4$. 
    \begin{itemize}
        \item For $n=2$, we get a lower bound only w.r.t. the quadratic module. We cannot hope for more, because the preordering is \emph{saturated} there, i.e., $\preo(\gcube) = \pos(\gcube)$ for $n=2$~\cite{scheiderer_sums_2006}.
        \item For $n=3$, $n=4$, we get lower bounds by taking a globally nonnegative (but not SOS) polynomial, and homogenizing  it. The strength of the bound is better if the degree of the polynomial is lower. Note that this follows the Hilbert cases for homogeneous polynomials: at degree $2$, all globally nonnegative forms are SOS. For $n=3$, this is true up to and including degree $4$.
        For $n \geq 4$, there are globally nonnegative forms of degree $4$ that are not SOS.
    \end{itemize}
    }
\end{remark}

\subsection{Transfer theorems, constraint qualifications and {\L}ojasiewicz exponents}

Given the convergence rate results on the hypercube, a natural question is whether such rates also hold for other ``simple'' sets such as the unit sphere, unit ball, and standard simplex. Rather than providing separate statements for each of these sets, we observe that all of their standard descriptions give rise to Archimedean quadratic modules, and satisfy the constraint qualification LICQ. The following~\Cref{thm:LICQupper} shows that this is sufficient to obtain an $O(k^{-4})$ rate of convergence. 
This theorem thus implies a rate of $O(k^{-4})$ for standard descriptions of the unit ball, the unit sphere and the simplex. 
The proof of this theorem relies crucially on the $O(k^{-4})$ convergence rate of the Putinar-based hierarchy on the hypercube stated in~\Cref{THM:cubeupper}.~\Cref{thm:LICQupper} can be viewed as transferring that result to any Archimedean quadratic module for which LICQ holds. 
We will provide a mathematical exposition of the proof of this result in a forthcoming version of this preprint. 
For completeness, note that the instance used in the lower degree bound for general Archimedean quadratic modules (\Cref{THM:generalputinarlower}) does not satisfy LICQ.

\begin{theorem}[{\lean}] \label{thm:LICQupper}
Let $\mg$ be a set of constraints defining an Archimedean quadratic module. Suppose that LICQ holds for $\mg$ (\Cref{def::LICQ}), and that $\se(\mg)$ is nonempty. Let $p \in \pos(\mg)$.
Then, there exist constants $C, k_0, \e_0 > 0$, depending only on $p$ and $\mg$, such that
\[
    \e_{p, k, \qm(\mg)} \leq C \cdot \frac{1}{k^4} \quad (\forall k \geq k_0); \quad 
    k_{p, \e, \qm(\mg)} \leq C \cdot \frac{1}{\sqrt[4]{\e}} \quad (\forall \e \leq \e_0).
\]

\end{theorem}

The convergence rate of the SOS hierarchy for a constraint set $\mg$ has been connected to the {\L}ojasiewicz exponent $\Lo(\mg)$ of $\mg$ at several places in the literature~\cite{Schweighofer2004,baldi2023effective,baldi2025lojasiewicz,heijmans2026degree}. It is known that if $\mg$  satisfies LICQ, we have $\Lo(\mg) = 1$~\cite{baldi2023effective}. It is thus a reasonable question whether the quartic convergence established above might extend to all (Archimedean) $\mg$ with ${\Lo(\mg) = 1}$. This turns out to be false, as the following example shows.

\begin{definition} \label{DEF:Loja}
    Let $\mg$ be a set of polynomial constraints defining a compact, nonempty semialgebraic set $\se(\mg)$. We write $\Lo(\mg)$ for the {\L}ojasiewicz exponent of $\mg$, which is the smallest $\Lo \geq 1$ so that the following holds: For every closed, bounded neighborhood $\mathcal{N} \supseteq \se(\mg)$, there exists a $C_\mathcal{N} > 0$ with
    \[
        \mathrm{dist}(\x, \se(\mg))^\Lo \leq C_\mathcal{N} \cdot \max\{0,\max_{g \in \mg} -g(\x)\} \quad \forall \x \in \mathcal{N}. 
    \]
\end{definition}

\begin{proposition}[{see, e.g.,~\cite{baldi2025lojasiewicz}}] 
    Let $\mg$ be a set of polynomial constraints. Assume that $\qm(\mg)$ is Archimedean and that $\mg$ satisfies LICQ (\Cref{def::LICQ}). Then $\Lo(\mg) = 1$.
\end{proposition}

\begin{example} \label{EXMP:MFCQ}
Let $n=3$ and write
\begin{equation*}
\begin{aligned}
q_1(x_1,x_2)&=x_1^2-3x_2^2,
&q_2(x_1,x_2)&=x_2^2-3x_1^2,\\
q_3(x_1,x_2)&=4x_1x_2-x_1^2-x_2^2,
&q_4(x_1,x_2)&=-4x_1x_2-x_1^2-x_2^2.
\end{aligned}
\end{equation*}
Consider the constraint set $\mg = (g_i : 1 \leq i \leq 5)$ given by
\[
g_i=x_3-q_i\quad(1\le i\le4),\quad
g_5=1-x_1^2-x_2^2-x_3^2.
\]
The quadratic module $\qm(\mg)$ is Archimedean. The {\L}ojasiewicz exponent of $\mg$ is $\Lo(\mg) = 1$. The LICQ for $\mg$ fails at $\x = (0,0,0)$. We have $x_3 \geq 0$ for all $\x \in \se(\mg)$. 
\end{example}

\begin{theorem}[{\lean}]
Let $n=3$ and consider the constraint set $\mg$ of~\Cref{EXMP:MFCQ}. Let $p = x_3$. Then $p \in \pos(\mg)$. There exist constants $c, k_0, \e_0 > 0$ such that,
\[
\e_{p, k, \qm(\mg)} \geq c \cdot \frac{1}{k^2} \quad (\forall k \geq k_0); 
\quad k_{p, \e, \qm(\mg)} \geq c \cdot  \frac{1}{\sqrt{\e}} \quad (\forall \e \leq \e_0).
\]
\end{theorem}

\begin{remark}
By combining \Cref{THM:cubeupper} with the techniques presented in~\cite{heijmans2026degree}, it follows that for any set of constraints $\mg$ defining an Archimedean quadratic module $\qm(\mg)$ and a nonempty semialgebraic set $\se(\mg)$, and any $p \in \pos(\mg)$, there exist constants $C, k_0, \e_0 > 0$, depending only on $p$ and $\mg$, such that
\[
    \e_{p, k, \qm(\mg)} \leq C \cdot \frac{1}{k^{2/\Lo(\mg)}} \quad (\forall k \geq k_0); \quad 
    k_{p, \e, \qm(\mg)} \leq C \cdot \frac{1}{\e^{\Lo(\mg)/2}} \quad (\forall \e \leq \e_0).
\]
In particular, a quadratic convergence-rate upper bound is obtained when $\Lo(\mg)=1$.
\end{remark}

\section{Historical development of convergence rates and degree bounds}
\label{sec:tables}
The following tables attempt to exhaustively survey the existing literature on \Cref{question}.
Each rate displayed below corresponds to $\e_{p,k,\mathcal C}$ as defined in~\Cref{eq:epkC}, with the cone $\mathcal C$ specified in each row. 
In the first four tables, the ``Rate'' column displays upper bounds as
$k\to\infty$, and, unless a row states otherwise, objectives are arbitrary polynomials; the stated bound holds up to constants which may depend on $\mg$ and $p$ (as noted earlier, many of the cited works make these constants more explicit). 
In the last table, the ``Rate'' column displays either lower bounds with the ``$\Omega$'' notation or exact error bounds. 
In~\Cref{tab:sos-box-rates}, hypercube rows use the standard quadratic generators, or the corresponding standard interval description under a bounded degree conversion. 
In~\Cref{tab:sos-special-rates}, for the ball, one has $\qm=\preo$ because there is only one inequality generator. 
For any compact univariate basic closed set, one also has $\qm=\preo$. 
In~\Cref{tab:sos-construction-rates}, the sparse rows assume variable cliques satisfying the so-called \emph{running intersection property}; the objective is a sum of polynomials supported on individual cliques, and each constraint and certificate term is supported on one such clique. 

In~\Cref{tab:sos-lower-bounds}, the univariate examples use their displayed nonstandard generators.
Their lower bounds do not apply to the interval described by
$1-x^2$. The Baldi--Slot objective~\cite[Thm.~4, Cor.~16]{BaldiSlot2024} is already certified in
$\preo(\gcube)_2$. The final row concerns SOS approximation and is not a lower bound
for the hypercube hierarchy.

\begingroup
\small
\setlength{\tabcolsep}{5pt}
\setlength{\extrarowheight}{3pt}
\renewcommand{\arraystretch}{1.18}
\setlength{\SOSwidth}{\dimexpr\linewidth-5\tabcolsep\relax}
\setlength{\LTcapwidth}{\linewidth}
\setlength{\LTpre}{12pt}
\setlength{\LTpost}{8pt}
\hypersetup{hidelinks}

\if{
\noindent
\textbf{Conventions.}
$E_r=f_{\min}-f_r$, with certificate degree at most $2r$.
$\qm$ denotes the Putinar quadratic module, $\preo$ the Schm\"udgen
preordering, and $\Sigma+\ideal$ sums of squares modulo the equality ideal.
The dimension, objective, and defining polynomials are fixed.
Unless specified otherwise, objectives are arbitrary polynomials.
A superscript $\mathrm{p}$ marks a reference cited as a preprint.
The first four tables give upper bounds; the fifth gives lower
bounds or exact errors.

\medskip\noindent
For general-set estimates, normalize the \L{}ojasiewicz exponent by
\[
 \operatorname{dist}(x,K)^L
 \le C\,\operatorname{viol}_{g,h}(x),
 \qquad L\ge1,
\]
on a fixed surrounding compact set. The violation function measures
negative inequality residuals and absolute equality residuals.
The exponent depends on the defining polynomials, not only on $K$.
Every $\qm$ upper bound also holds for $\preo$, since $\qm_{2r}\subseteq\preo_{2r}$.
}\fi


\begin{longtable}{
    Z{.21\SOSwidth}
    Z{.16\SOSwidth}
    Z{.18\SOSwidth}
    Z{.2\SOSwidth}
    Z{.2\SOSwidth}
}
\caption{General compact and Archimedean sets}
\label{tab:sos-general-rates}\\
\toprule
\rowcolor{white}
\textbf{Setting}
 & \textbf{Certificate}
 & \textbf{Rate $O(\cdot)$}
 & \textbf{Note}
 & \textbf{Reference} \\
\midrule
\endfirsthead

\caption[]{General compact and Archimedean sets (continued)}\\
\toprule
\rowcolor{white}
\textbf{Setting}
 & \textbf{Certificate}
 & \textbf{Rate $O(\cdot)$}
 & \textbf{Note}
 & \textbf{Reference} \\
\midrule
\endhead

\midrule
\multicolumn{5}{r}{\footnotesize Continued on the next page.} \\
\endfoot

\bottomrule
\endlastfoot

\rowcolor{SOSstripe}
Compact $\se(\mg)$
 & $\preo(\mg)$
 & $k^{-\alpha_{\mg}}$
 & $\alpha_{\mg}>0$
 & \cite[Thms.~3--4]{Schweighofer2004} \\

\rowcolor{white}
Archimedean $\qm(\mg)$
 & $\qm(\mg)$
 & $(\log k)^{-\alpha_{\mg}}$ 
 & $\alpha_{\mg}>0$
 & \cite[Thm.~6]{NieSchweighofer2007} \\

\rowcolor{SOSstripe}
Archimedean $\qm(\mg)$
 & $\qm(\mg)$
 & $k^{-2/(5n \Lo(\mg))}$
 & $n\ge2$
 & \cite[Thm.~1.7, Eq.~(21)]{baldi2023effective} \\

\rowcolor{white}
Archimedean $\qm(\mg)$

 & $\qm(\mg)$
 & $k^{-1/(7\Lo(\mg)+3)}$
 & 
 & \cite[Cor.~3.3]{baldi2025lojasiewicz} \\

\rowcolor{SOSstripe}
Archimedean $\qm(\mg)$

 & $\qm(\mg)$
 & $k^{-1/(2\Lo(\mg))}$
 & 
 & \cite[Thm.~4]{heijmans2026degree} \\

\rowcolor{white}
Compact $\se(\mg)$
 & $\preo(\mg)$
 & $k^{-1/\Lo(\mg)}$
 & 
 & \cite[Thm.~6]{heijmans2026degree} \\

\rowcolor{SOSstripe}
Compact $\se(\mg)$
 & $\preo(\mg)$
 & $k^{-1/\Lo(\mg)}$
 & with normalization, via a subcone
 & \cite[Thm.~4.7, Cor.~4.8]{tran2026convergence} \\

\rowcolor{white}
Compact $\se(\mg)$; with explicit hypercube generators
 & $\qm(\mg)$
 & $\dfrac{(\log k)^{3/(2\Lo(\mg))}}{k^{1/\Lo(\mg)}}$
 & 
 & \cite[Thm.~3.7]{tran2026extension} \\

\end{longtable}





\begin{longtable}{
    Z{.21\SOSwidth}
    Z{.16\SOSwidth}
    Z{.18\SOSwidth}
    Z{.2\SOSwidth}
    Z{.2\SOSwidth}
}
\caption{Hypercube}
\label{tab:sos-box-rates}\\
\toprule
\rowcolor{white}
\textbf{Setting}
 & \textbf{Certificate}
 & \textbf{Rate $O(\cdot)$}
 & \textbf{Note}
 & \textbf{Reference} \\
\midrule
\endfirsthead

\caption[]{Hypercube: convergence-rate upper bounds (continued)}\\
\toprule
\rowcolor{white}
\textbf{Setting}
 & \textbf{Certificate}
 & \textbf{Rate $O(\cdot)$}
 & \textbf{Note}
 & \textbf{Reference} \\
\midrule
\endhead

\midrule
\multicolumn{5}{r}{\footnotesize Continued on the next page.} \\
\endfoot

\bottomrule
\endlastfoot

\rowcolor{SOSstripe}
Hypercube
 & $\preo(\gcube)$ 
 & $k^{-1}$
 & via Handelman representation
 & \cite[Thm.~1.4]{deKlerk2010error} \\

\rowcolor{white}
Hypercube
 & $\qm(\gcube)$
 & $(\log k)^{-1}$
 & \emph{conditional} on a conjecture of de~Klerk--Laurent.
 & \cite[Thm.~4]{magron2015error} \\

\rowcolor{SOSstripe}
Hypercube
 & $\preo(\gcube)$
 & $k^{-2}$
 & 
 & 
 \cite[Thm.~2]{LaurentSlot2023} \\

\rowcolor{white}
Hypercube
 & $\qm(\gcube)$
 & $k^{-1}$
 & 
 & 
        \cite[Thm.~3, Cor.~15]{BaldiSlot2024}
    \\

\rowcolor{SOSstripe}
Hypercube
 & $\qm(\gcube)$
 & $(\log k)/k$
 & SOS approximation
 & \cite[Cor.~24]{gribling2026revisiting}
    \\

\rowcolor{white}
Hypercube
 & $\qm(\gcube)$
 & $(\log k)^3/k^2$
 & SOS approximation
 & \cite[Thm.~7]{gribling2026} \\

\rowcolor{SOSstripe}
Hypercube
 & $\preo(\gcube)$ 
 & $k^{-2}$
 & trigonometric SOS
 & \cite[Thm.~3.1]{bach2023exponential} \\
\end{longtable}



\begin{longtable}{
    Z{.21\SOSwidth}
    Z{.16\SOSwidth}
    Z{.18\SOSwidth}
    Z{.2\SOSwidth}
    Z{.2\SOSwidth}
}
\caption{Sphere, ball, simplex, polytopes, and univariate sets}
\label{tab:sos-special-rates}\\
\toprule
\rowcolor{white}
\textbf{Setting}
 & \textbf{Certificate}
 & \textbf{Rate $O(\cdot)$}
 & \textbf{Note}
 & \textbf{Reference} \\
\midrule
\endfirsthead

\caption[]{Sphere, ball, simplex, polytopes, and univariate sets
(continued)}\\
\toprule
\rowcolor{white}
\textbf{Setting}
 & \textbf{Certificate}
 & \textbf{Rate $O(\cdot)$}
 & \textbf{Note}
 & \textbf{Reference} \\
\midrule
\endhead

\midrule
\multicolumn{5}{r}{\footnotesize Continued on the next page.} \\
\endfoot

\bottomrule
\endlastfoot

\rowcolor{SOSstripe}
Sphere
 & $\qm(\mg_{\rm sph})$
 & $k^{-1}$
 & forms of even degree
 & \cite[Thm.~3.12]{reznick1995uniform} \\

\rowcolor{white}
Sphere
 & $\qm(\mg_{\rm sph})$
 & $k^{-1}$
 & forms of even degree
 & \cite[Cor.~7.2]{doherty2012convergence} \\ 

\rowcolor{SOSstripe}
Sphere
 & $\qm(\mg_{\rm sph})$
 & $k^{-2}$
 & forms of even degree $2d$;
original theorem assumes $d\le n$
 & \cite[Thm.~2]{fang2021sum} \\

\rowcolor{white}
Sphere
 & $\qm(\mg_{\rm sph})$
 & $k^{-2}$
 & forms of even degree; via moment-relaxation
 & \cite[Thms.~3.1--3.2]{blomenhofer2026moment} \\ 

\rowcolor{SOSstripe}
Sphere
 & $\qm(\mg_{\rm sph})$
 & $k^{-2}$
 & extends~\cite{fang2021sum} to arbitrary polynomials
 & \cite[Thm.~4.15, Cor.~4.16]{tran2026convergence}\\

\rowcolor{white}
Unit ball
 & $\qm=\preo$
 & $k^{-2}$
 & 
 & \href{https://arxiv.org/abs/2111.04610}
        {\cite[Thm.~3]{slot2022sum}} \\

\rowcolor{SOSstripe}
Simplex
 & $\preo$ 
 & $k^{-1}$
 & P\'olya-type representation
 & \cite[Thm.~1]{powers2001new};  \cite[Thm.~1.1]{de2006ptas} \\

\rowcolor{white}
Simplex
 & LP-based approximation
 & $k^{-1}$
 & quadratic objectives
 & \cite[Thm.~3.2, Cor.~3.1]{bomze2002solving}\\

\rowcolor{SOSstripe}
Simplex
 & $\preo$
 & $k^{-2}$
 & 
 & \cite[Thm.~4]{slot2022sum} \\

\rowcolor{white}
Compact polytope
 & $\preo$
 & $k^{-1}$
 & linear inequality description
 & \cite[Thm.~4.14]{tran2026convergence}\\

\rowcolor{SOSstripe}
Compact univariate basic closed set
 & $\qm = \preo$
 & $k^{-2}$
 & arbitrary given generators
 & \cite[Thm.~1]{henrion2026convergence} \\


\end{longtable}



\begin{longtable}{
    Z{.21\SOSwidth}
    Z{.16\SOSwidth}
    Z{.18\SOSwidth}
    Z{.2\SOSwidth}
    Z{.2\SOSwidth}
}
\caption{Products, liftings, affine transfers, and sparse constructions}
\label{tab:sos-construction-rates}\\
\toprule
\rowcolor{white}
\textbf{Setting}
 & \textbf{Certificate}
 & \textbf{Rate $O(\cdot)$}
 & \textbf{Note}
 & \textbf{Reference} \\
\midrule
\endfirsthead

\caption[]{Products, lifting, affine transfers, and sparse constructions
(continued)}\\
\toprule
\rowcolor{white}
\textbf{Setting}
 & \textbf{Certificate}
 & \textbf{Rate $O(\cdot)$}
 & \textbf{Note}
 & \textbf{Reference} \\
\midrule
\endhead

\midrule
\multicolumn{5}{r}{\footnotesize Continued on the next page.} \\
\endfoot

\bottomrule
\endlastfoot

\rowcolor{SOSstripe}
Product of two spheres
 & $\preo$
 & $k^{-2}$
 & 
 & \cite[Thm.~3]{magron2026convergence}\\

\rowcolor{white}
Any fixed finite product of spheres
 & $\preo$
 & $k^{-2}$
 & 
 & \cite[Thm.~10]{magron2026convergence}\\

\rowcolor{SOSstripe}
Products of balls, simplices, hypercubes, and compatible mixtures
 & $\preo$
 & $k^{-2}$
 & 
 & \cite[Cor.~12]{magron2026convergence}; \cite[Thm.~3.3]{tran2026convergence}\\

\rowcolor{white}
Two spheres
 & Matrix SOS, hence $\preo$
 & $k^{-2}$
 & $f(x,y)=y^{\mathsf T}F(x)y$, with $F$ a form of even degree
 & \cite[Thm.~2]{fang2021sum}\\

\rowcolor{SOSstripe}
Algebraic sections of products
 & $\preo$
 & $k^{-1/\Lo(\mg)}$
 & with suitable quadratic-rate kernels
 & \cite[Thm.~4.3]{tran2026convergence}\\

\rowcolor{white}
Tensor products of compatible positive polynomial kernels
 & $\preo$
 & Factor errors combine; slowest factor rate is retained
 & 
 & \cite[Thm.~11]{magron2026convergence}\\

\rowcolor{SOSstripe}
Positive lifting into a hypercube
 & Original $\qm$ or $\preo$
 & Transfers hypercube degree bounds with controlled degree growth
 & followed by polynomial substitution
 & \cite[Thm.~2]{heijmans2026degree}\\

\rowcolor{white}
Hypercube
 & Sparse $\preo(\gcube)$
 & $k^{-2/(\omega+3)}$
 & $\omega$ is the largest variable-clique size
 & \cite[Thm.~6]{korda2025convergence}\\

\rowcolor{SOSstripe}
Block-Archimedean sparse problem
 & Sparse $\qm$
 & $k^{-\beta_{\rm sp}}$;
   if all $\Lo_j=1$:
   $k^{-18/(55\omega+238)}$
 & $\beta_{\rm sp}$ depends on clique sizes and $\Lo_j$'s (local Łojasiewicz exponents)
 & \cite[Thm.~8]{korda2025convergence}\\  

\rowcolor{white}
Invertible affine images
 & Transported $\qm$ or $\preo$
 & Same exponent
 & ellipsoids, parallelepipeds, affine simplices
 & Elementary degree-preserving affine substitution \\

\end{longtable}




\begin{longtable}{
    Z{.21\SOSwidth}
    Z{.16\SOSwidth}
    Z{.18\SOSwidth}
    Z{.2\SOSwidth}
    Z{.2\SOSwidth}
}
\caption{Lower bounds and exact-error constructions}
\label{tab:sos-lower-bounds}\\
\toprule
\rowcolor{white}
\textbf{Setting}
 & \textbf{Certificate}
 & \textbf{Rate $\Omega(\cdot)$}
 & \textbf{Note}
 & \textbf{Reference} \\
\midrule
\endfirsthead

\caption[]{Lower bounds and exact-error constructions (continued)}\\
\toprule
\rowcolor{white}
\textbf{Setting}
 & \textbf{Certificate}
 & \textbf{Rate $\Omega(\cdot)$}
 & \textbf{Note}
 & \textbf{Reference} \\
\midrule
\endhead

\midrule
\multicolumn{5}{r}{\footnotesize Continued on the next page.} \\
\endfoot

\bottomrule
\endlastfoot

\rowcolor{SOSstripe}
${g=(1-x^2)^3}$; ${p = 1-x^2}$
 & $\qm(g)=\preo(g)$
 & $k^{-2}$
 & 
 & \cite[Thm.~4]{Stengle1996} \\

\rowcolor{white}
${g=(1-x^2)^3}$; ${p = 1-x^2}$
 & $\qm(g)=\preo(g)$
 & $\dfrac{1}{k(k-2)}$
 &  $k\ge3$, exact rate
 & \cite[Thm.~4]{Henrion2026} \\

\rowcolor{SOSstripe}
$g=(1-x)^3(1+x)$; $p=1-x$
 & $\qm(g)=\preo(g)$
 & $\dfrac{1}{2k(k-1)}$
 & $k\ge2$, exact rate  & \cite[Thm.~2]{henrion2026convergence}\\

\rowcolor{white}
$\mg=\{x_1^3-x_2^2,\;1-x_1\}$; $p=x_1$
 & $\qm(\mg)$
 & $\dfrac{1}{(2k-1)(2k-3)}$ 
 & $k\ge2$, exact rate
 & \cite[Rmk.~3]{henrion2026convergence}\\

\rowcolor{SOSstripe}
Standard hypercube; $p=(1-x_1^2)(1-x_2^2)$
 & $\qm(\gcube)$
 & $k^{-8}$
 & $n\ge2$;
 & \cite[Thm.~4, Cor.~16]{BaldiSlot2024}\\

\rowcolor{white}
Standard hypercube
 & $\Sigma$ 
& $k^{-2}$ 
& $n \ge 2$: SOS approximation
 & \cite[Thm.~10]{gribling2026revisiting}
\end{longtable}

\endgroup


\begin{thebibliography}{GdKV26b}

\bibitem[Art27]{artin_uber_1927}
Emil Artin.
\newblock Uber die {Zerlegung} definiter {Funktionen} in {Quadrate}.
\newblock {\em Abhandlungen aus dem Mathematischen Seminar der Universitat
  Hamburg}, 5(1):100--115, December 1927.

\bibitem[BdK02]{bomze2002solving}
Immanuel~M Bomze and Etienne de~Klerk.
\newblock Solving standard quadratic optimization problems via linear,
  semidefinite and copositive programming.
\newblock {\em Journal of Global Optimization}, 24(2):163--185, 2002.

\bibitem[BK26]{BaldiKocvara:PC}
Lorenzo Baldi and Michal Kočvara.
\newblock Conjectures in real algebra and polynomial optimization through high
  precision semidefinite programming, 2026.
\newblock In preparation.

\bibitem[BL26]{blomenhofer2026moment}
Alexander~Taveira Blomenhofer and Monique Laurent.
\newblock {Moment-sos and spectral hierarchies for polynomial optimization on
  the sphere and quantum de Finetti theorems}.
\newblock {\em SIAM Journal on Optimization}, 36(1):204--232, 2026.

\bibitem[BM23]{baldi2023effective}
Lorenzo Baldi and Bernard Mourrain.
\newblock {On the effective Putinar’s Positivstellensatz and moment
  approximation}.
\newblock {\em Mathematical Programming}, 200(1):71--103, 2023.

\bibitem[BMP25]{baldi2025lojasiewicz}
Lorenzo Baldi, Bernard Mourrain, and Adam Parusi{\'n}ski.
\newblock {On {\L}ojasiewicz inequalities and the effective Putinar's
  Positivstellensatz}.
\newblock {\em Journal of Algebra}, 662:741--767, 2025.

\bibitem[BR23]{bach2023exponential}
Francis Bach and Alessandro Rudi.
\newblock Exponential convergence of sum-of-squares hierarchies for
  trigonometric polynomials.
\newblock {\em SIAM Journal on Optimization}, 33(3):2137--2159, 2023.

\bibitem[BS24]{BaldiSlot2024}
Lorenzo Baldi and Lucas Slot.
\newblock Degree bounds for {P}utinar's {P}ositivstellensatz on the hypercube.
\newblock {\em SIAM J. Appl. Algebra Geom.}, 8(1):1--25, 2024.

\bibitem[dKL10]{deKlerk2010error}
Etienne de~Klerk and Monique Laurent.
\newblock Error bounds for some semidefinite programming approaches to
  polynomial minimization on the hypercube.
\newblock {\em SIAM Journal on Optimization}, 20(6):3104--3120, 2010.

\bibitem[dKLP06]{de2006ptas}
Etienne de~Klerk, Monique Laurent, and Pablo~A Parrilo.
\newblock {A PTAS for the minimization of polynomials of fixed degree over the
  simplex}.
\newblock {\em Theoretical Computer Science}, 361(2-3):210--225, 2006.

\bibitem[DW12]{doherty2012convergence}
Andrew~C Doherty and Stephanie Wehner.
\newblock {Convergence of SDP hierarchies for polynomial optimization on the
  hypersphere}.
\newblock {\em arXiv preprint arXiv:1210.5048}, 2012.

\bibitem[FF21]{fang2021sum}
Kun Fang and Hamza Fawzi.
\newblock The sum-of-squares hierarchy on the sphere and applications in
  quantum information theory.
\newblock {\em Mathematical Programming}, 190(1):331--360, 2021.

\bibitem[GdKV26a]{gribling2026revisiting}
Sander Gribling, Etienne de~Klerk, and Juan Vera.
\newblock {Revisiting the convergence rate of the Lasserre hierarchy for
  polynomial optimization over the hypercube}.
\newblock {\em Optimization}, pages 1--28, 2026.

\bibitem[GdKV26b]{gribling2026}
Sander Gribling, Etienne de~Klerk, and Juan~C. Vera.
\newblock {Squared polynomial approximation kernels for the hypercube: improved
  error bounds and implications for Lasserre hierarchies}.
\newblock {\em arXiv preprint}, 2026.

\bibitem[Hen26]{Henrion2026}
Didier Henrion.
\newblock Solving {Stengle's} example in rational arithmetic: Exact values of
  the {Moment-SOS} relaxations.
\newblock {\em Open Journal of Mathematical Optimization}, 7:1--18, 2026.

\bibitem[HKVZ26]{heijmans2026degree}
Olga Heijmans-Kuryatnikova, Juan~C Vera, and Luis~F Zuluaga.
\newblock {Degree Bounds for Positivstellensätze of general semialgebraic
  sets}.
\newblock {\em arXiv preprint arXiv:2605.15821}, 2026.

\bibitem[HSED26]{henrion2026convergence}
Didier Henrion and Mohab Safey El~Din.
\newblock Convergence rate of the moment-sos hierarchy for univariate
  polynomial optimization.
\newblock {\em arXiv preprint arXiv:2609.20544}, 2026.

\bibitem[KMRZ25]{korda2025convergence}
Milan Korda, Victor Magron, and Rodolfo Rios-Zertuche.
\newblock Convergence rates for sums-of-squares hierarchies with correlative
  sparsity.
\newblock {\em Mathematical Programming}, 209(1):435--473, 2025.

\bibitem[Kri64]{krivine_anneaux_1964}
J.~L. Krivine.
\newblock Anneaux préordonnés.
\newblock {\em Journal d’Analyse Mathématique}, 12(1):307--326, December
  1964.

\bibitem[Las01]{lasserre2001global}
Jean~B Lasserre.
\newblock Global optimization with polynomials and the problem of moments.
\newblock {\em SIAM Journal on optimization}, 11(3):796--817, 2001.

\bibitem[LS23]{LaurentSlot2023}
Monique Laurent and Lucas Slot.
\newblock An effective version of {S}chm\"udgen's {P}ositivstellensatz for the
  hypercube.
\newblock {\em Optim. Lett.}, 17(3):515--530, 2023.

\bibitem[Mag15]{magron2015error}
Victor Magron.
\newblock {Error bounds for polynomial optimization over the hypercube using
  Putinar type representations}.
\newblock {\em Optimization Letters}, 9(5):887--895, 2015.

\bibitem[Mag26]{magron2026convergence}
Victor Magron.
\newblock Convergence rates for polynomial optimization on set products.
\newblock {\em SIAM Journal on Optimization}, 36(3):1867--1886, 2026.

\bibitem[NS07]{NieSchweighofer2007}
Jiawang Nie and Markus Schweighofer.
\newblock On the complexity of {Putinar}'s {Positivstellensatz}.
\newblock {\em Journal of Complexity}, 23(1):135--150, 2007.

\bibitem[Par00]{parrilo2000thesis}
Pablo~A Parrilo.
\newblock {\em Structured semidefinite programs and semialgebraic geometry
  methods in robustness and optimization}.
\newblock California Institute of Technology, 2000.

\bibitem[PR01]{powers2001new}
Victoria Powers and Bruce Reznick.
\newblock {A new bound for P{\'o}lya's theorem with applications to polynomials
  positive on polyhedra}.
\newblock {\em Journal of pure and applied algebra}, 164(1-2):221--229, 2001.

\bibitem[Put93]{Putinar1993}
Mihai Putinar.
\newblock Positive polynomials on compact semi-algebraic sets.
\newblock {\em Indiana Univ. Math. J.}, 42(3):969--984, 1993.

\bibitem[Rez95]{reznick1995uniform}
Bruce Reznick.
\newblock {Uniform denominators in Hilbert's seventeenth problem}.
\newblock {\em Mathematische Zeitschrift}, 220(1):75--97, 1995.

\bibitem[Rez07]{Reznick2007}
Bruce Reznick.
\newblock On {H}ilbert's construction of positive polynomials.
\newblock {\em Preprint}, 2007.

\bibitem[Sch91]{Schmudgen1991}
Konrad Schm{\"u}dgen.
\newblock The {$K$}-moment problem for compact semi-algebraic sets.
\newblock {\em Math. Ann.}, 289:203--206, 1991.

\bibitem[Sch04]{Schweighofer2004}
Markus Schweighofer.
\newblock On the complexity of {Schm{\"u}dgen}'s {Positivstellensatz}.
\newblock {\em Journal of Complexity}, 20(4):529--543, 2004.

\bibitem[Sch06]{scheiderer_sums_2006}
Claus Scheiderer.
\newblock Sums of squares on real algebraic surfaces.
\newblock {\em Manuscripta Mathematica}, 119(4):395--410, April 2006.

\bibitem[Slo22]{slot2022sum}
Lucas Slot.
\newblock {Sum-of-squares hierarchies for polynomial optimization and the
  Christoffel--Darboux kernel}.
\newblock {\em SIAM Journal on Optimization}, 32(4):2612--2635, 2022.

\bibitem[Ste74]{stengle_nullstellensatz_1974}
G.~Stengle.
\newblock A nullstellensatz and a positivstellensatz in semialgebraic geometry.
\newblock {\em Mathematische Annalen}, 207(2):87--97, June 1974.

\bibitem[Ste96]{Stengle1996}
Gilbert Stengle.
\newblock Complexity estimates for the {Schm{\"u}dgen} {Positivstellensatz}.
\newblock {\em J. Complexity}, 12(2):167--174, 1996.

\bibitem[TT26a]{tran2026extension}
Hoang~Anh Tran and Kim-Chuan Toh.
\newblock {Extension on the convergence rates of moment-SOS hierarchies via
  approximation of truncated moment sequences}.
\newblock {\em arXiv preprint arXiv:2609.37286}, 2026.

\bibitem[TT26b]{tran2026convergence}
Hoang~Anh Tran and Kim-Chuan Toh.
\newblock {On the convergence rates of moment-SOS hierarchies approximation of
  truncated moment sequences}.
\newblock {\em Mathematical Programming}, pages 1--61, 2026.

\end{thebibliography}

\end{document}